\documentstyle[12pt]{article}
\begin{document}

\textwidth18cm
\textheight21.5cm
%\topmargin0.5cm
%\oddsidemargin0.5cm
%\evensidemargin1cm
%\pagestyle{empty}
%\textwidth18cm
%\textheight28cm
%\countdef\revised=100
%\revised=1

\newtheorem{Ex}{\bf Example}%[section]
\newtheorem{Rem}{\bf Remark}%[section]
\newtheorem{Sa}{\bf Proposition}%[section]
\newtheorem{Le}[Sa]{\bf Lemma}
\newtheorem{Hi}[Sa]{\bf Lemma}
\newtheorem{Th}[Sa]{\bf Theorem}
\newtheorem{Fo}[Sa]{\bf  Corollary}
\newtheorem{Beh}[Sa]{\bf Proposition}
\newtheorem{De}{\bf Definition}%[section]

\newcommand{\DS}{ \mathrel{\mathop{+}\limits^.   }}
\newcommand{\RR}{{\rm I\kern-0.14em R }}
\newcommand{\NN}{{\rm l\kern-0.14em N}}
\newcommand{\CC}{{\rm\raise 0.192ex\vbox{\hrule height
1.22ex width 0.8pt}\kern-0.29em C}}
\newcommand{\ZZ}{ {\sf Z}\hspace{-0.4em}{\sf Z}\ }
\newcommand{\TT}{ { /}\hspace{-0.4em}{|} \ }
\newcommand{\AAAAA}{ \hspace{-0.6em}{\bf a}\ }
\def\ii{{\rm i}}
\def\dd{{\, \rm d }}
%\begin{flushleft}
%\end{flushleft}

%\begin{flushleft}
%e-mail-address \\
%X.400:    leitenberger@mathematik.uni-leipzig.d400.de
%\end{flushleft}

\begin{center}
{\LARGE
Quantum Lobachevsky Planes}
\end{center}

\begin{flushleft}
{\large FRANK LEITENBERGER } \\
{\small \it Fachbereich Mathematik, Universit\"at Rostock,
Rostock, D-18051, Germany. \\
e-mail: frank.leitenberger@mathematik.uni-rostock.de }
\end{flushleft}

%\begin{flushleft}
%
%\end{flushleft}

%\renewcommand{\baselinestretch}{2}
%\small \normalsize

\vspace{0.5cm}
{\small
{\bf Abstract:} We classify all $SL(2,\RR)$-covariant
Poisson structures on the
Lo\-ba\-chev\-sky plane with respect to all
multiplicative Poisson structures on $SL(2,\RR)$ and describe
Quantisations for all these Poisson structures.}
\vspace{0.5cm}

02.20.Qs; 03.65.Fd

\newpage
\noindent
{\Large \bf I. Introduction}
\setcounter{equation}{0}
\vspace{0.5cm}

\noindent
A first step into the direction of a quantisation of
the notions
of algebraic geometry (cf. Ref. 1) is the quantisation
of the three fundamental one-dimensional complex
domains:
the complex plane, the Riemannian sphere
and the complex upper half plane.
Quantum planes were considered in Ref. 2.
Quantum Riemannian spheres are discussed in
Ref. 3 and Ref. 4.
A quantum upper half plane appears firstly in Ref. 5
(for certain other quantisations of the upper half plane
see also Remark 7).
In this letter we consider the case of the quantum
complex upper half plane.

In the theory of quantum groups
the classification of certain Poisson structures gives
generally a
good insight into the problem of
the classification of quantum
structures.
At first we give the full solution of
the classification problem in the
classical limit (i.e. the description of all possible
$SL(2,\RR)$-covariant
Poisson structures for the upper half plane with respect to
the action of all possible
multiplicative Poisson structures on $SL(2,\RR)$).
Further we describe classes of quantum structures, which
reproduce all listed Poisson structures in the classical
limit.
We obtain two-parameter quantisations
of the upper half plane
for every action of one of the quantum groups $SL_q (2,\RR )$,
$SU_q (1,1 )$ and $SL_h (2,\RR )$.
\vspace{0.5cm}

\noindent
{\Large \bf II. Preliminaries}
\vspace{0.5cm}

\noindent
In this section we introduce the basic concepts
(see also Ref. 5, Ref. 6 and Ref. 7).
Let $({\rm A},\bigtriangleup , \epsilon)$
be a Hopf algebra. A \it left
quantum space \rm
$({\rm H},\phi )$ is an algebra H along with
an algebra homomorphism
$\phi: {\rm H} \rightarrow {\rm A} \otimes
{\rm H} $
such that
$ (\bigtriangleup \otimes id) \phi = (id \otimes \phi ) \phi $
and
$ (\epsilon \otimes id) \phi = id  $.
Two left quantum spaces $({\rm H}_1,\phi_1 ),\ ({\rm H}_2,\phi_2 )$
are \it isomorphic \rm if
there is an algebra isomorphism $h:{\rm H}_1
\rightarrow {\rm H}_2$, such that
$\phi_2 \circ h = (id \otimes h) \circ \phi_1$.
Further consider commutative Hopf algebras
and commutative quantum spaces.
We say that the Poisson bracket
$\{.,.\}$ on A is \it multiplicative \rm if
$ \{ \bigtriangleup (x),
\bigtriangleup (y) \}_{{\rm A}\otimes {\rm A}} =
  \bigtriangleup (\{ x,y \}) $, where
  ${\rm A}\otimes {\rm A}$ carries the
Poisson structure of the direct product.
We call a Poisson bracket
$\{.,.\}$ on H  \it covariant \rm if
$ \{ \phi (x), \phi (y) \}_{{\rm A}\otimes {\rm H}} =
  \phi (\{ x,y \}) $, where ${\rm A}\otimes {\rm H}$
  carries the direct
product structure from A and H.
We say that two multiplicative or covariant, respectively,
Poisson brackets on A and, respectively, H  are
{\it equivalent} if they intertwine
with an automorphism
of A and, respectively, H.
\vspace{0.5cm}

\noindent
{\Large \bf  III. The classical case}
\setcounter{equation}{0}
\vspace{0.5cm}

\noindent
{\bf A. Poisson-Structures on $SL(2,\RR)$}

\noindent
Let  ${\rm A}
=\RR [a,b,c,d]/(ad-bc-1)$ be the commutative
unital    algebra of
polynomial functions in the coordinates of
\[  SL(2,\RR) =
\left\{ \left. \left(
\begin{array}{cc}
a & b \\ c & d
\end{array}
\right) \ \right| \ \  a,b,c,d \in \RR,\ ad-bd=1
\right\} .\]
${\rm A}$ becomes an Hopf algebra with respect to the
Hopf multiplication
\[ \bigtriangleup
\left(
\begin{array}{ll} a & b \\ c & d
\end{array}
\right)
=
\left(
\begin{array}{ll}
a \otimes a + b \otimes c  &
a \otimes b + b \otimes d  \\
c \otimes a + d \otimes c  &
c \otimes b + d \otimes d
\end{array}
\right)
\]
and the counit $\epsilon$ with
\[ \epsilon
\left(
\begin{array}{ll} a & b \\ c & d
\end{array}
\right)
=
\left(
\begin{array}{ll}
1  &
0  \\
0  &
1
\end{array}
\right).   \]
A has the coinverse $\kappa$ with
\[ \kappa
\left(
\begin{array}{ll} a & b \\ c & d
\end{array}
\right)
=
\left(
\begin{array}{cc}
d  &
- b  \\
- c &
a
\end{array}
\right).
\]

We consider three types of
multiplicative Poisson algebras
${\rm A}^A_\lambda,$
${\rm A}^K_\lambda $ and
${\rm A}^N$:

\noindent
1. $ {\rm A}_\lambda^A= ({\rm A}, \{, \}_\lambda) $,
$\lambda \in \RR,\ \lambda\neq 0$, with
\[   \begin{array}{llc}
 \{a,b\}_\lambda & = &  \lambda ab,  \\
 \{a,c\}_\lambda & = &  \lambda ac,  \\
 \{a,d\}_\lambda & = &  2 \lambda bc, \\
 \{b,c\}_\lambda & = &   0, \\
 \{b,d\}_\lambda & = &  \lambda bd,  \\
 \{c,d\}_\lambda & = &  \lambda cd;
\end{array}
 \]
2. $ {\rm A}_\lambda^K = ({\rm A}, \{, \}_\lambda) $,
$ \lambda > 0 $, with
\[   \begin{array}{llc}
 \{a,b\}_\lambda & = &  \lambda (1 - a^2-b^2),  \\
 \{a,c\}_\lambda & = &  \lambda (a^2 + c^2-1),  \\
 \{a,d\}_\lambda & = &  \lambda (a-d)(b-c), \\
 \{b,c\}_\lambda & = &  \lambda (a+d)(b+c), \\
 \{b,d\}_\lambda & = &  \lambda (b^2+d^2-1),  \\
 \{c,d\}_\lambda & = &  \lambda (1-c^2-d^2)
\end{array}
 \]
and \\
3. $ {\rm A}^N = ({\rm A}, \{, \}) $
with
\[   \begin{array}{llc}
 \{a,b\} & = &   (1 - a^2),  \\
 \{a,c\} & = &   c^2,        \\
 \{a,d\} & = &   c (d-a),    \\
 \{b,c\} & = &   c (d+a),    \\
 \{b,d\} & = &   (d^2 - 1),  \\
 \{c,d\} & = & - c^2.
\end{array}
 \]

$A_\lambda^A$ ($\lambda\in\RR,\ \lambda\neq 0$),
     $A_\lambda^K$ ($\lambda > 0$) and
     $A^N$
we shortly denote by
     $A_\lambda^I$ ($I=A,K,N$).

\begin{Rem} \rm The notation ${\rm A}^I_\lambda$ with
$I=A,K,N$ is
justified by the fact, that
the Poisson brackets
of
${\rm A}^A_\lambda,
 {\rm A}^K_\lambda$ and
${\rm A}^N$, respectively,
vanish on the
subgroups of the {\rm KAN}-decomposition
\[G_A=
\left\{ \left. \left(
\begin{array}{ll} e^x & 0 \\ 0 & e^{-x} \end{array}
\right) \ \right| \ \ x \in \RR    \right\} ,
\]
\[G_K=
\left\{ \left. \left( \begin{array}{rl} \cos x & \sin x \\
- \sin x & \cos x \end{array}
\right) \ \right| \ \ x \in [0,2\pi )
 \right\}, \] and
\[G_N=
\left\{ \left. \left(
\begin{array}{ll} 1 & x \\ 0 & 1 \end{array}
\right) \ \right| \ \ x \in \RR    \right\}
\]
respectively.
\end{Rem}

\begin{Sa}  Every nontrivial multiplicative Poisson
structure on {\rm A} is equivalent
to one of the structures
${\rm A}_\lambda^A\  (\lambda \neq 0),
\ {\rm A}_\lambda^K\ (\lambda > 0)$,
${\rm A}^N$.
All these structures are
non-equivalent.
\end{Sa}

{\it Proof.} We give a scetch proof.
Set $g=sl(2,\RR)$. Because
$H^1(g,g\wedge g)=0$, every
Poisson structure on $SL(2,\RR)$ arises from a classical
r-matrix $r\in g\wedge g$, which satisfies the
modified classical
Yang Baxter equation
$(\wedge_1^3 ad_x) C(r)=0, \ \forall x \in g$
with a certain element $C(r)\in g\wedge g\wedge g$
(cf. Ref. 6).
Because $g=sl(2,\RR)$ is a simple Lie algebra,
the one-dimensional representation $\wedge^3_1 ad$ is
trivial. That is, the MCYBE is satisfied for every $r\in g$.

To show, that two structures are equivalent, it is enough to show
that their r-matrices are connected by an automorphism of
%the Hopf algebra $({\rm A},\bigtriangleup,\epsilon )$.
$sl(2,\RR)$.
All automorphisms are generated by inner automorphisms
(i.e. the adjoint action of $SL(2,\RR)$) and the automorphism
$\alpha$ with
\[ \alpha \left( \begin{array}{cc} a & b  \\  c & d \end{array}
\right) = \left( \begin{array}{cc} a & -b \\ -c & d \end{array}
\right).
\]

Let
$
e_{-1} = \left( \begin{array}{cc}  0 & 1 \\ 0 & 0  \end{array} \right),
e_{0}  = \left( \begin{array}{cc}  1 & 0 \\ 0 & -1 \end{array} \right),
e_{1}  = \left( \begin{array}{cc}  0 & 0 \\ 1 & 0  \end{array} \right)$
and  \\
$r=  \alpha e_1     \wedge e_{-1}
       + \beta  e_0 \wedge e_1
       + \gamma e_0 \wedge e_{-1}$.
Studying the adjoint action $\wedge_1^2Ad$ of $SL(2,\RR)$
and the action of $\alpha$ on
$r\in g\wedge g$ we recognize,
that $r$ is equivalent to one of the
following elements
\begin{eqnarray*}
r^A_\lambda = & \lambda  e_{1} \wedge e_{-1}, &
\lambda \in \RR,\ \lambda\neq 0 \\
r^K_\lambda = &
\lambda (e_0 \wedge e_{1} + e_0 \wedge e_{-1}), & \lambda > 0, \\
r^N = &  e_0 \wedge e_1, &
\end{eqnarray*}
All these elements are non-equivalent.

Because A consists of the matrix elements of the finite
dimensional representations $\rho:SL(2, \RR) \rightarrow \RR^n$,
one can give the Poisson structures by the formula
\[
\{
l_1(\rho_1 (g) v_1),
l_2(\rho_2 (g) v_2)
\} =
(l_1\otimes l_2) ([ (\rho_1 \otimes \rho_2) (r_0),
\rho_1(g) \otimes \rho_2(g) ] v_1 \otimes v_2)
\]
($v \in \RR^d, l\in (\RR^{d})^*$, cf. Ref. 8).
The structures
$ {\rm A}^A_{\lambda},
\ {\rm A}^K_{\lambda},
\ {\rm A}^N,   $
correspond to the above normal forms
$r^A_\lambda,\ r^K_\lambda, r^N$ for $r$.
$\Box$
\vspace{0.5cm}

\noindent
{\bf B. Poisson structures on the upper half plane}

\noindent
Further consider the subalgebra
${\rm H} \subseteq {\rm A}$ generated by the elements
 $A= a^2 + b^2,\ B= ac+bd, \ D= c^2+d^2$
(cf. Ref. 5) and the left coaction
$\phi = \bigtriangleup|_{\rm H}: {\rm H} \rightarrow
{\rm A} \otimes {\rm H} $ with

\[  \begin{array}{lll}
 \phi (A) & = & a^2 \otimes A + b^2 \otimes D + 2ab \otimes B,  \\
 \phi (D) & = & c^2 \otimes A + d^2 \otimes D + 2cd \otimes B , \\
 \phi (B) & = & ac  \otimes A + bd  \otimes D + (ad+bc) \otimes B.
\end{array}
 \]

By the next proposition we give a
complete classification of all possible
covariant Poisson structures on H with respect
to ${\rm A}_\lambda^I$.

\begin{Sa}
\rm 1. \it  All  with respect to
${\rm A}^I_\lambda$
covariant Poisson structures
are given  by
the one-parameter series:

(i) $ {\rm H}^A_{\lambda,\mu }$
   ($  \lambda\in\RR $, $\lambda\neq 0$, $\mu \in\RR$)

\begin{eqnarray*} \{ A,B \} & = & 2 \lambda A (B+\mu), \\
                  \{ A,D \} & = & 4 \lambda B (B+\mu), \\
                  \{ B,D \} & = & 2 \lambda D (B+\mu),
\end{eqnarray*}

(ii) ${\rm H}^K_{\lambda,\mu }$
    ($\lambda >0$, $\mu\in \RR$)

\begin{eqnarray*} \{ A,B \} & = & 2 \lambda A (A+D+\mu), \\
                       \{ A,D \} & = & 4 \lambda B (A+D+\mu), \\
                       \{ B,D \} & = & 2 \lambda D (A+D+\mu),
\end{eqnarray*}

(iii) ${\rm H}^N_\mu$
     ($\mu\in\RR$)

\begin{eqnarray*}  \{ A,B \} & = & 2  A (D+\mu), \\
                         \{ A,D \} & = & 4  B (D+\mu), \\
                         \{ B,D \} & = & 2  D (D+\mu)
\end{eqnarray*}
for $I=n$.

{\rm 2}. Two Poisson algebras
${\rm H}_{\lambda ,\mu_1}^I$ and
${\rm H}_{\lambda ,\mu_2}^I$
are equivalent if and only if $|\mu_1|=|\mu_2|$.
\end{Sa}

\noindent
{\it Proof.} 1. All left invariant Poisson stuctures
(i.e. $\phi (\{ x , y \} )
=  \{\phi (x) ,\phi (y) \}$,
$\forall x,y\in {\rm H}$,
where ${\rm A} \otimes {\rm H}$ carries a Poisson structure
which is the direct product of the zero structure on A and of the
structure on H)
are given by
\begin{eqnarray*}   \{ A,B \} &   =   &   2 \mu A,   \\
                    \{ A,D \} &   =   &   4 \mu B,   \\
                    \{ B,D \} &   =   &   2 \mu D,
\end{eqnarray*}
$(\mu\in\RR)$.
We obtain the above formulas by the
calculation for a fixed $\mu$
(for example $\mu=0$)
and from the fact that the difference
of two left covariant structures is
left invariant (cf. Ref. 7).

2. The proof of the equivalence follows from the fact that
$({\rm H},\phi )$ has the unique automorphism
$x\rightarrow -x$, $\forall x\in {\rm H}$.
$\Box$

The following Proposition shows, that we can realize the
Poisson algebras ${\rm H}_{\lambda,\mu}^I$ as subalgebras
of ${\rm A}_\lambda^I$.

Consider subalgebras
${\rm H}_{\alpha,\beta,\gamma} \subset {\rm A} $
generated by
elements $\overline{A},\overline{B},\overline{D} \in {\rm A}$ with
\begin{eqnarray*}
\overline{A}: & = & \alpha a^2  +  \beta b^2  + 2\gamma ab, \\
\overline{B}: & = & \alpha ac   +  \beta bd   + \gamma (ad+bc),  \\
\overline{D}: & = & \alpha c^2  +  \beta d^2  + 2\gamma cd.
\end{eqnarray*}
We have ${\rm H}_{1,1,0} = {\rm H}$ and
the correspondence between $A,B,D$ and the overlined elements,
arranges an isomorphism
${\rm H}_{\alpha,\beta,\gamma} \cong {\rm H}$
(i.e. $\overline{A}\, \overline{D}-\overline{B}^2=1$),
if and only if $\alpha\beta-\gamma^2=1$.

\begin{Sa}
Let $\alpha \beta - \gamma^2 = 1 $. Then
${\rm H}_{\alpha, \beta, \gamma}$ is a Poisson
subalgebra
of ${\rm A}^I_{\lambda}$ and we have

\[ \begin{array}{ll}
(i) \ \ \ \  & {\rm H}_{\alpha,\beta,\gamma}
\cong {\rm H}^A_{\lambda,-\gamma } , \\
(ii)\ \ \ \  & {\rm H}_{\alpha,\beta,\gamma}
\cong
{\rm H}^K_{\lambda, -\alpha-\beta}
\ \
and \\
(iii) \ \ \ \  & {\rm H}_{\alpha,\beta,\gamma}
\cong {\rm H}^N_{ \lambda , -\beta}.
\end{array}
\]
\end{Sa}

\noindent
{\it Proof}. The isomorphies can be
verified by an explicite calculation.

\begin{Rem} \rm The
preceeding Proposition admits the realisation of all
${\rm H}_{\lambda,\mu}^I$
by ${\rm H}_{\alpha,\beta,\gamma}$
with real $\alpha,\beta,\gamma$ except the cases
${\rm H}^N_{0}$ and
${\rm H}^K_{\lambda,\mu}$, $\vert \mu \vert <2$.
For example
${\rm H}^A_{\lambda,\mu} \cong {\rm H}_{1,1+\mu^2,\mu}$,
${\rm H}^K_{\lambda,\mu} \cong
{\rm H}_{ \frac{\mu + \sqrt{\mu^2 - 4}}{2},
    \frac{\mu - \sqrt{\mu^2 - 4}}{2}, 0}$ and
${\rm H}^N_{\mu}
\cong {\rm H}_{-\frac{1}{\mu},-\mu,0}$.
In the other cases
%${\rm H}_{\lambda,0}^n$
we have only complex realisations, for example
${\rm H}^K_{\lambda,\mu} \cong
{\rm H}_{ \frac{\mu +\ii \sqrt{4- \mu^2 }}{2},
    \frac{\mu -\ii \sqrt{4- \mu^2 }}{2}, 0}$ and
$
{\rm H}^N_{0}
\cong
{\rm H}_{0,0,\ii}
$. The realisations are not unique.

%\begin{eqnarray*}
%\overline{A}: & = & 2\ii  ab,      \\
%\overline{B}: & = &  \ii (ad+bc), \\
%\overline{D}: & = & 2\ii  cd.
%\end{eqnarray*}
\end{Rem}
\vspace{0.5cm}
\noindent
\bf C. Geometric interpretation \rm

\noindent
There is an interpretation of the parameter $\mu$ in terms of
fixed points $z_0$ of the upper half plane.

Let $z \rightarrow \frac{az+b}{cz+d}$ be the
usual action of $SL(2,\RR)$ on
the complex coordinate $z=x+iy$, $y>0$
of the upper half plane, i.e.
\begin{eqnarray*}
 x & \rightarrow & \frac{ ac(x^2+y^2) + (ad+bc)x + bd }{
                       c^2(x^2+y^2) +   2cd x  + d^2}  , \\
 y & \rightarrow & \frac{y}{
                        c^2(x^2+y^2) +   2cd x  + d^2}.
\end{eqnarray*}
Further let $G_{z_0}$ be the subgroup of $SL(2,\RR)$, which fixes
the point $z_0=x_0+\ii y_0$, $y_0 > 0$. It follows that
the functions
\[ x=x(a,b,c,d)=
   \frac{ ac(x_0^2+y_0^2) + (ad+bc)x_0 + bd }{
                       c^2(x_0^2+y_0^2) +   2cd x_0  + d^2}  \]
and
\[   y=y(a,b,c,d)=
   \frac{y_0}{ c^2(x_0^2+y_0^2) +   2cd x_0  + d^2}   \]
are left $G_{z_0}$-invariant and can be identified with functions on
$G_{z_0}\setminus G$. We recognize that
$\overline{B} = x y^{-1}$,
$\overline{D} =   y^{-1}$ generate
${\rm H}_{ y_0 + \frac{x_0^2}{y_0}, \frac{x_0}{y_0}, \frac{1}{y_0} }$.
For the parameter $\mu$ we get the interpretation
\begin{eqnarray*}
\mu =  & {\rm ctg} ( {\rm arg} ({\it z_0}) )
              = \frac{ {\rm Re } z_0 }{ {\rm Im} z_0}, &
{\rm if} \ \  I=A \\
\mu =  &  - \frac{1+\mid z_0 \mid^2}{ {\rm Im} z_0},  &
{\rm if} \ \  I=K \\
\mu =  &  - \frac{1}{ {\rm Im} z_0 }, &
{\rm if} \ \  I=N.
\end{eqnarray*}

We remark, that $x,y   \notin {\rm H}  =\RR
[\overline{A},\overline{B},\overline{D} ]$ and
                $\overline{A},\overline{B},\overline{D}
\notin {\rm H}':=\RR [x,y]  $, but
H and ${\rm H}'$
are dense subalgebras of the algebra of smooth functions
on the upper half plane (H and H$'$ are dense on every compact
subset).
With respect to the xy-coordinates the Poisson brackets have the
form
\[  \{x,y\} = -2\lambda (x  y         + \mu y^2)  \]
if $I=A$,
\[  \{x,y\} = -2\lambda (x^2y+y+y^3   + \mu y^2)  \]
if $I=K$ and
\[  \{x,y\} = -2        (   y         + \mu y^2)  \]
if $I=N$.
\vspace{0.5cm}

\noindent
{ \Large \bf  IV. The quantum case }
\vspace{0.5cm}
\setcounter{equation}{0}

\noindent
{\bf A. The Quantum groups ${\cal A}_h^I$}

\noindent
First we recall the definition of the three known quantum
deformations of $SL(2,\RR)$.

The Poisson algebras ${\rm A}_\lambda^I$
correspond to  the quantum algebras
${\cal A}_h^A$  $(h\in \RR,\ h\neq 2\pi n,\ n\in\NN)$,
${\cal A}_h^K$  $(h>0)$ and
$  {\cal A}^N$. For all these algebras
we write shortly ${\cal A}_h^I$ ($I=A,K,N$).
The algebras ${\cal A}_h^I$
are generated by elements $a,b,c,d$
and relations
\[ \begin{array}{rll}
   ab & = & e^{\ii h} ba , \\
   ac & = & e^{\ii h} ca , \\
   bd & = & e^{\ii h} db , \\
   cd & = & e^{\ii h} dc , \\
   bc & = & cb,     \\
   ad-da & = & (e^{\ii h} - e^{-\ii h}) bc,     \\
   ad-e^{\ii h} bc & = & da- e^{-\ii h} cb=1, \ \  (h\neq 2\pi n),
\end{array}  \] for $I=A$,

\[ \begin{array}{rll}
   ab-ba & = &  \ \ii h   (1-a^2-b^2) , \\
   ac-ca & = &  \ \ii h   (a^2+c^2-1) , \\
   bd-db & = &  \ \ii h   (b^2+d^2-1), \\
   cd-dc & = &  \ \ii h   (1-c^2-d^2), \\
   ad-da & = &  \ \ii h   (ab+ba+ac+ca+bd+db+cd+dc), \\
   bc-cb & = &  \ \ii h   (ab+ba-ac-ca-bd-db+cd+dc), \\
   ad+da-bc-cb & = &    \ 2+ h (2-a^2-b^2-c^2-d^2),\ \ \ (h>0),
\end{array}  \]
for $I=K$ and

\[ \begin{array}{rcl}
   ab-ba & = & \ii   (1-a^2), \\
   ac-ca & = & \ii   c^2, \\
   bd-db & = & \ii   (d^2-1), \\
   cd-dc & = & -\ii  c^2, \\
   ad-da & = & \ii   c(d-a), \\
   bc-cb & = & \ii   (dc+ca), \\
   ad-bc+\ii ac & = &  da-cb-\ii ca=1.
\end{array}  \] for $I=N$.

The algebras ${\cal A}^I_h$
become  Hopf algebras with respect to the
Hopf multiplication
\[ \bigtriangleup
\left(
\begin{array}{ll} a & b \\ c & d
\end{array}
\right)
=
\left(
\begin{array}{ll}
a \otimes a + b \otimes c  &
a \otimes b + b \otimes d  \\
c \otimes a + d \otimes c  &
c \otimes b + d \otimes d
\end{array}
\right)
\]
and Hopf $*$-algebras with respect to the involution \\
$a^*:=a,\ b^*:=b,\ c^*:=c,\ d^*:=d$.

The counit $\epsilon$ is given by

\[ \epsilon
\left(
\begin{array}{ll} a & b \\ c & d
\end{array}
\right)
=
\left(
\begin{array}{ll}
1  &
0  \\
0  &
1
\end{array}
\right)
\]
and the coinverse $\kappa$ is given by
\[ \kappa
\left(
\begin{array}{ll} a & b \\ c & d
\end{array}
\right)
=
\left(
\begin{array}{cc}
d  &
-q^{-1} b  \\
-q c &
a
\end{array}
\right)
\] for $I=A$,

\begin{eqnarray*}
\kappa(a) & = &   \frac{1}{4} \left((2- \frac{2-2h^2}{1+h^2})a
                  +(2+ \frac{2-2h^2}{1+h^2})d
                  -\frac{2h}{1+h^2}b
                  -\frac{2h}{1+h^2}c \right),              \\
\kappa(b) & = &   \frac{1}{4}\left(\frac{2}{1+h^2} a
                  - \frac{2}{1+h^2} d
                  -(2+\frac{2-2h^2}{1+h^2})b
                  +(2-\frac{2-2h^2}{1+h^2})c \right),         \\
\kappa(c) & = &  \frac{1}{4}\left(\frac{2}{1+h^2} a
                  -  \frac{2}{1+h^2} d
                  +  (2- \frac{2-2h^2}{1+h^2})b
                  -  (2+ \frac{2-2h^2}{1+h^2})c   \right),     \\
\kappa(d) & = &  \frac{1}{4}\left((2+ \frac{2-2h^2}{1+h^2})a
                  +(2- \frac{2-2h^2}{1+h^2})d
                  +\frac{2h}{1+h^2}b
                  +\frac{2h}{1+h^2}c\right)
\end{eqnarray*}
for $I=K$ and

\[ \kappa
\left(
\begin{array}{ll} a & b \\ c & d
\end{array}
\right)
=
\left(
\begin{array}{cc}
d+  c                &
-b+(d-a+ c)  \\
-c  &
a-c
\end{array}
\right)
\] for $I=N$.

\begin{Rem} \rm ${\cal A}^K_h$ is equivalent to the quantum group
              $SU_q (1,1)$, $q\in \RR$
which is given by elements $\alpha,\alpha^*, \beta,\beta^*$
and relations

\[ \begin{array}{rcl}
   \alpha\beta  & = & q \beta \alpha , \\
   \alpha\beta^* & = & q \beta^*\alpha , \\
   \beta\beta^* & = & \beta^*\beta , \\
   \beta\alpha^* & = & q \alpha^*\beta , \\
   \beta^*\alpha^* & = & q \alpha^*\beta^* , \\
   \alpha  \alpha^* -q^{2} \beta\beta^* & = &
   \alpha^*\alpha -  \beta^*\beta =1.
\end{array}  \]
The formulas of $SU_q (1,1)$ and ${\cal A}^K_h$ are connected by
the transformation of the deformation parameter
$ h= \frac{1-q}{1+q} $ and by the
"quantum Cayley transformation" of the matrix elements
 \begin{eqnarray*}
 \alpha \   & := &  \frac{1}{2}  (a+d +\ii (b-c)),  \\
 \alpha^*  & := &  \frac{1}{2}   (a+d +\ii (-b+c)), \\
 \beta \    & := &  \frac{1}{2}  (a-d -\ii (b+c)),  \\
 \beta^*   & := &  \frac{1}{2}   (a-d +\ii (b+c)).
\end{eqnarray*}
\end{Rem}
\begin{Rem} \rm
Instead of ${\cal A}^N$ one  considers Hopf $*$-algebras
$SL_h(2,\RR),\\ (h\neq 0)$ where the
algebra structure is replaced by

$  ab-ba  =  \ii h  (1-a^2),
   ac-ca  =  \ii  h c^2,
   bd-db  =  \ii  h (d^2-1),$

$  cd-dc  =  -\ii h c^2,
   ad-da  =  \ii  h c(d-a),
   bc-cb  =  \ii  h (dc+ca), $

$  ad-bc+\ii h ac  =   da-cb-\ii h ca=1   $
%\end{array}  \]
\\ (cf. Ref. 9 and Ref. 10).
 All these structures are equivalent to ${\cal A}^N$ by the
Hopf $*$-algebra isomorphism
$ a\rightarrow a,\
  b\rightarrow hb,\
  c\rightarrow \frac{1}{h} c,\
  d\rightarrow d $.
\end{Rem}

\vspace{0.5cm}
\noindent
{\bf B. The quantum spaces $({\cal H}_{h,k,d}^I,\phi)$ }

\noindent
{\bf 1. The algebras ${\cal H}_{h,k,d}^I$ }

\noindent
The Poisson algebras ${\rm H}_{\lambda,\mu}^I$ correspond
to the algebras
${\cal H}_{h,k,d}^a$ $(h\neq 2\pi n)$,
${\cal H}_{h,k,d}^k$ $(h>0)$,
${\cal H}_{k,d}^n$,
$k,d \in \RR$.
For all these algebras we write shortly
${\cal H}_{h,k,d}^I$.

We define the
${\cal H}_{h,k,d}^I$ as algebras
which are generated by elements $A,$$B,$$C,$$D$
and relations
\[ \begin{array}{ccl}
C  & = & e^{-\ii h}  B + k(1-e^{-\ii h})  , \\
AB & = & e^{2\ii h} BA + k(1-e^{2\ii h}) A , \\
%AD - DA  & = & (e^{2\ii h} - e^{-2\ii h})
%BC  ???????? +\frac{ k (1+e^{2h}) }{e^{2h} } B, \\
BD & = & e^{2\ii h} DB + k(1-e^{2\ii h}) D , \\
AD - e^{2\ii h}  BC & = & d + k(1-e^{2\ii h} ) B, \\
DA - e^{-2\ii h} CB & = & d + k(1-e^{-2\ii h}) B,\ \ \
(h\neq 2\pi n),
\end{array} \]
for $I=A$,

\[ \begin{array}{cll}
C & = & B - \ii h (A + D + k ) , \\
AB-BA  & = & \ii h (2kA+2A^2 +AD +BC) , \\
% AD-DA  & = & \ii h (2k(B+C)+AB+BA+AC+CA+2BD+2DC) , \\
BD-DB  & = & \ii h (2kD+2D^2 +AD +CB), \\
AD-B^2  & = & d + \ii h (kB+AB+BD), \\
DA-C^2  & = & d - \ii h (kC+CA+DC),\ \ \ (h>0),
\end{array} \]
for $I=K$ and

\[ \begin{array}{cll}
C      & = & B - \ii (D+k), \\
AB-BA  & = & 2\ii  A(D+k) - 2 B(D+k) , \\
% AD-DA  & = & 4\ii  B(D+k) - 2 (3D^2 + 4k D +k^2 ) , \\
BD-DB  & = & 2\ii  D(D+k), \\
AD-BC  & = & d + 2\ii  B(D+k)  , \\
DA-CB  & = & d - 2\ii  (D+k)C
\end{array} \]
for $I=N$.

The algebras
${\cal H}_{h,k,d}^I$ become $*$-algebras
with respect to the involutions
\[  A^* := A, \ \ D^* := D, \ \
B^* := C = e^{-\ii h} B + k (1-e^{-\ii h}) ,\ \ C^*:=B   \]
for $I=A$,
\[  A^* := A, \ \ D^* := D, \ \  B^* := C= B - \ii h (A+D +k),\ \
C^* := B  \]
for $I=K$ and
\[  A^* := A, \ \ D^* := D, \ \  B^* := C = B - \ii  (D + k),\ \
C^* := B  \]
for $I=N$.

\begin{Rem} \rm Formally we recover the Poisson structures
of subsection III.B in the limit $t\rightarrow\infty$,
if we set $B=C$, $q=\lambda t$, $k=\mu $, $d=1$ and
$\{x,y\}:= \lim \frac{1}{it} [x,y]$.
\end{Rem}

Next we show, that we can realize the $*$-algebras
${\cal H}_{h,k,d}^I$
as $*$-subalgebras of ${\cal A}_h^I$. Let

\begin{eqnarray*}
\overline{A}: & = & \alpha a^2  +  \beta b^2  + \gamma (ab+ba), \\
\overline{B}: & = & \alpha ac   +  \beta bd   + \gamma (ad+bc), \\
\overline{C}: & = & \alpha ca   +  \beta db   + \gamma (da+cb), \\
\overline{D}: & = & \alpha c^2  +  \beta d^2  + \gamma (cd+dc).
\end{eqnarray*}

\begin{Sa} Let $\alpha,\beta,\gamma\in \RR$.
$\overline{A}, \overline{B}, \overline{C}, \overline{D}$
satisfy the relations  of

(i)
%$H_{\alpha,\beta,\gamma} \cong
$ {\cal H}^A_{h,\gamma,\alpha\beta -\gamma^2 }$,

(ii)
%$H_{\alpha,\beta,\gamma} \cong
${\cal H}^K_{h,-\alpha -\beta, \alpha\beta -\gamma^2 }$
and

(iii)
%$H_{\alpha,\beta,\gamma} \cong
${\cal H}^N_{h,-\beta,\alpha\beta -\gamma^2 }$.
\end{Sa}

\noindent
{\it Proof}.  The proof can be given by an explicite calculation.

\begin{Rem} \rm
Proposition 4 admits the realisation of all
${\cal H}^I_{h,k,d}$,
with real $\alpha,\beta,\gamma$ without the cases
${\cal H}^K_{h,0,d}$,  $|k| < 2$ and
${\cal H}^N_{h,0,d}$, $     d > 0$.
For example
     ${\cal H}^A_{h,k,d} \cong {\cal H}_{1,d+k^2,k}$,
${\cal H}^K_{h,k,d} \cong
{\cal H}_{ \frac{k + \sqrt{k^2 - 4}}{2},
    \frac{k - \sqrt{k^2 - 4}}{2}, 0}$ and
${\cal H}^N_{h,k,d} \cong {\cal H}_{-\frac{d}{k},-k,0}$.
In other cases
%${\rm H}_{h,0,d}^n$, $d>0$
we have only complex realisations, for example \\
${\rm H}^K_{h,k,d} \cong
H_{ \frac{k +\ii \sqrt{4- k^2 }}{2},
    \frac{k -\ii \sqrt{4- k^2 }}{2}, 0}$ and
${\cal H}_{0,0,\ii\sqrt{d} } \cong
{\cal H}^N_{h, 0, d}$.
%\begin{eqnarray*}
%\overline{A}: & = & 2\ii  ab,      \\
%\overline{B}: & = &  \ii (ad+bc), \\
%\overline{D}: & = & 2\ii  cd.
%\end{eqnarray*}
\end{Rem}

\noindent
\begin{Rem} \rm
1. The special case ${\cal H}_{1,1,0} \cong {\cal H}^A_{h,0,1}$
was first mentioned in Ref. 5, p. 188.
\\ 2. The formulas for
${\cal H}^A_{h,k,d}$ are similar those
from  Podles' sphere  $C(X_{\mu, \lambda, \rho})$ (cf. Ref. 4).
Both can be specified from the complex quantum space
${\Xi}_{q,\lambda,\rho}$ of $SL_q(2,\CC)$ from Ref. 11
by fixing certain involutions.
\\ 3. In  Ref. 12
the quantum spaces  ${\cal H}_{\alpha,\beta,0}
\cong {\cal H}^A_{h,0,\alpha\beta}$, $\alpha,\beta \geq 0$
were considered.
\\ 4. The quantum spaces ${\cal H}^K_{h,k,d}$ correspond
to one-parameter series of quantum discs $C_{\mu,q} (\overline{U})$
in Ref. 13. Formally the correspondence is
arranged by the quantum Cayley transformation (cf. Remark 3).
\end{Rem}

\vspace{0.5cm}
\noindent
{\bf 2. The coaction of ${\cal A}_h^I$ on
${\cal H}_{h,k,d}^I$ }

\noindent
We will describe an coaction $\phi$, i.e. an homomorphism
$\phi: {\cal H}^I_{h,k,d} \rightarrow {\cal A}_h^I \otimes
{\cal H}^I_{h,k,d}$
with
$ (\bigtriangleup \otimes id) \phi = (id \otimes \phi ) \phi $
and
$ (\epsilon \otimes id) \phi = id  $
such that we can call $({\cal H}_{h,k,d}^I,\phi) $
a quantum space with respect to ${\cal A}_h^I$.

According to Remark 6
consider ${\cal H}^I_{h,k,d}$ as a subalgebra of ${\cal A}^I_h$,
i.e. we identify $A,B,C,D$ with
$
\overline{A},
\overline{B},
\overline{C},
\overline{D}
$.
We obtain (independent from $I,h,k,d$)
\[ \begin{array}{lll}
\bigtriangleup(A)
& = & a^2 \otimes A+b^2 \otimes D+ab\otimes B+ba\otimes C, \\
\bigtriangleup(D)
& = & c^2 \otimes A+d^2 \otimes D+cd\otimes B+dc\otimes C, \\
\bigtriangleup(B)
& = & ac  \otimes A+bd  \otimes D+ad\otimes B+bc\otimes C, \\
\bigtriangleup(C)
& = & ca  \otimes A+db  \otimes D+cb\otimes B+da\otimes C.
   \end{array}
\]
That is, ${\cal H}^I_{h,k,d}$ is a left coideal of
${\cal A}^I_h$
and we have proven the following proposition.

\begin{Sa} The homomorphism
$\phi: {\cal H}^I_{h,k,d} \rightarrow {\cal A}^I_h
\otimes {\cal H}^I_{h,k,d}:
\phi:=\bigtriangleup|_{{\cal H}^I_{h,k,d}} $ defines
a left coaction, that is $({\cal H}_{h,k,d}^I,\phi)$
are left quantum spaces.
\end{Sa}

\begin{Rem} \rm
Every ${\cal H}_{h,k,d}^I$ is equivalent to one of the
quantum spaces
\begin{eqnarray*}
{\cal H}_{h,k,\pm 1}^I,
{\cal H}_{h,k,0}^I
\end{eqnarray*}
with  $k\geq 0$.
%$I=a,k,n$, $h \in [0,2\pi)$ if $I=a$,
%$h\in [0,2\pi)$ if $I=k,n$ and
%$ k \in [ 0,\infty )$.
\end{Rem}

\it Proof. \rm
We  achieve
$k\in [0,\infty )$
because of the automorphism $X\rightarrow -X$,
$X\in {\cal H}^I_{h,k,d}$ ,
and we achieve $d=0,\pm 1$
by the reparametrisation $X\rightarrow
\frac{1}{ \sqrt{|d|} }X$,
$X\in {\cal H}^I_{h,k,d}$.
$\Box$

%%%%   \[  xy= e^{2ih} yx + k'y^2 \]
%%%%
%%%%   \[          \]
%%%%
%%%%   \[ xy-yx= 2ihy +k'y^2   \]

% \begin{thebibliography}{777}
%{\bf References }

\vspace{0.5cm}
\small

% \bibitem[M]{M}
\noindent ${}^1 \rm Y$u.I. Manin,
{\it Quantum Groups and Non-commutative
Geometry}
(Centre de Re\-cher\-ches
Mathematiques, Montreal, 1988).

% \bibitem[S]{S}
\noindent ${}^{2}\rm K.$ Schm\"udgen, A. Sch\"uler,
%"Covariant
%differential calculi on quantum spaces and on quantum groups",
C.R. Acad. Sci. Paris {\bf 316}, Serie I, 1155 (1993).

% \bibitem[D]{D}
\noindent ${ }^3\rm M.S.$ Dijkhuizen, T.H. Koornwinder,
%"Quantum Homogeneous Spaces,
%Duality, and Quantum 2-Spheres",
Geom. Dedicata {\bf 52}, 291 (1994).

% \bibitem[P]{P}
\noindent ${}^{4}\rm P$. Podles,
%"Quantum Spheres",
{ Lett. Math. Phys.}
{\bf 14}, 193 (1987).

% \bibitem[F]{F}
\noindent ${}^5 \rm L.$ Faddeev, N. Reshetikhin, L. Takhtajan,
%"Quantization of Lie groups and Lie algebras",
{ Algebra i Analis } {\bf 1},1, 178 (1989).

% \bibitem[Dr]{Dr}
\noindent ${{ }^6 \rm V.}$G. Drinfeld
%{ "Quantum groups",}
(in  Proceedings of the
International Con\-gress of Mathematicians,
Berkeley, 1986).
%A.M. Gleason (ed.), 798-820,
%American Mathematical Society, Providence, RI

% \bibitem[L]{L}
\noindent ${}^7 \rm J.-H.$ Lu, A. Weinstein,
%"Poisson Lie groups,
%Dressing transformations and Bruhat decompositions",
J. Diff. Geom. {\bf 31}, 501 (1990).

%\bibitem[V]{V}
\noindent ${}^{8} \rm L$.L. Vaksman, J.S. Soibelman,
%"Algebra of
%functions on the quantum group $SU(2)$",
{ Funct. Anal. Appl.} {\bf 22},3, 170 (1988).

\noindent ${}^{9} \rm B.$A. Kupershmidt,
% "The quantum group $GL_h (2)$",
{ J. Phys. A: Math. Gen.} {\bf 25}, L1239 (1992).

%\bibitem[Z]{Z}
\noindent ${}^{10}\rm S.$ Zakrzewski,
%"A Hopf $*$-algebra
%of polynomials on the
%quantum $SL(2,\RR)$
%for a "{unitary} R-matrix",
Lett. Math. Phys. {\bf 22}, 287 (1991).

%%%%%        \bibitem[ ${\ \ }^1$]{A}
\noindent ${ }^{11} \rm J.$ Apel, K. Schmuedgen,
%"Clas\-si\-fi\-ca\-tion of
%Three-\-Di\-men\-si\-o\-nal
%Covariant Differential Calculi on Podles'
%Quantum Spheres and Related Spaces",
Lett. Math. Phys. {\bf 32}, 25 (1994).

% \bibitem[N]{N}
\noindent ${}^{12} \rm M.$ Noumi,
{\it Macdonald's symmetric polynomials as
zonal spherical functions on some homogeneous spaces},
{to appear in Adv. in Math.}.

% \bibitem[K]{K}
\noindent ${}^{13} \rm S.$ Klimek, A. Lesniewski,
%"A Two-Parameter
%Quantum Deformation of the Unit Disc",
{ J. Func. Anal.} {\bf 115}, 1 (1993).

% \bibitem[Ta]{Ta} Takhtajan, L.A.,
%{\it Introduction to Quantum groups},
%in {\it Quantum groups}, Proceedings of the
%8th Int. Workshop on Math.^M
%Phys. in Clausthal,^M
%ed. by Doebner, H.-D., Lect. Notes Phys. 370, 1990.^M

%\bibitem[C]{ ${\ \ }^{14}$ }
%\noindent ${ }^2\rm C$hari, V., Pressley, A.: ^M
%{\it A Guide to Quantum Groups,} ^M
%Cambridge University Press, Cambridge, 1994, p. 651.^M
% \bibitem[Co]{Co} Connes, A., Noncommutative^M
%Differential geometry, Academic Press 1994.^M
%\end{thebibliography}

\end{document}